\documentclass[a4paper,10pt,reqno]{amsart}
\usepackage{dsfont,amsmath,amsfonts,amscd,amssymb,graphicx,mathrsfs,eufrak,upgreek,pdflscape,colonequals,euscript,amsmath}
\usepackage[usenames]{color}
\usepackage{setspace}
\usepackage{array,multirow,booktabs,longtable}
\usepackage{comment}
\usepackage{cancel}
\usepackage{upref}
\renewcommand*\partial{\textnormal{\reflectbox{6}}}
\usepackage{nccmath}%%for medium sized math
\usepackage{calc}

\usepackage{array,multirow,booktabs,longtable}
\newcommand{\marginparstretch}{0.6}
\let\oldmarginpar\marginpar
\renewcommand\marginpar[1]{\-\oldmarginpar[\framebox{\setstretch{\marginparstretch}\begin{minipage}{\marginparwidth}{\raggedleft\tiny #1}\end{minipage}}]{\framebox{\setstretch{\marginparstretch}\begin{minipage}{\marginparwidth}{\raggedright\tiny #1}\end{minipage}}}}

\newcommand{\xdownarrow}[1]{%
  {\left\downarrow\vbox to #1{}\right.\kern-\nulldelimiterspace}
}

\usepackage[pagebackref]{hyperref}

\hypersetup{
  colorlinks   = true,          %Colors links instead of ugly boxes
  urlcolor     = blue,          %Color for external hyperlinks
  linkcolor    = purple,          %Color of internal links
  citecolor   = blue             %Color of citations
}

\newtheorem{thm}{Theorem}[section]
\newtheorem{theorem}[thm]{Theorem}

\newtheorem{lemma}[thm]{Lemma}

\theoremstyle{definition} 
\newtheorem{definition}[thm]{Definition}

\usepackage{amsopn}

\usepackage{upgreek}
\usepackage{colonequals}
\usepackage{tikz}
\usetikzlibrary{arrows,decorations.pathmorphing,decorations.pathreplacing,positioning,shapes.geometric,shapes.misc,decorations.markings,decorations.fractals,calc,patterns}

\tikzset{
        DB/.style={circle,draw=black,circle,fill=black,inner sep=0pt, minimum size=4pt},
        DW/.style={circle,draw=black,inner sep=0pt, minimum size=4pt},
        gap/.style={inner sep=0.5pt,fill=white}}

\newcommand{\Spec}{\mathop{\rm Spec}\nolimits}

\usepackage{euscript}
\newcommand{\scrA}{\EuScript{A}}
\newcommand{\scrB}{\EuScript{B}}

\newcommand{\scrJ}{\EuScript{J}}

\newcommand{\scrL}{\EuScript{L}}

\newcommand{\scrN}{\EuScript{N}}
\newcommand{\scrO}{\EuScript{O}}
\newcommand{\scrP}{\EuScript{P}}

\newcommand{\scrR}{\EuScript{R}}
\newcommand{\scrS}{\EuScript{S}}
\newcommand{\scrT}{\EuScript{T}}
\newcommand{\scrU}{\EuScript{U}}

\newcommand{\scrX}{\EuScript{X}}
\newcommand{\scrY}{\EuScript{Y}}
\newcommand{\scrZ}{\EuScript{Z}}

\newcommand{\Jac}{\scrJ\mathrm{ac}}
\newcommand\Curve{\mathtt{C}}
\newcommand{\lcl}{(\kern -2.5pt(}
\newcommand{\rcl}{)\kern -2.5pt)}
\newcommand{\llangle}{\langle\kern -2.5pt\langle}
\newcommand{\rrangle}{\rangle\kern -2.5pt\rangle}
\newcommand{\ringtwo}{\mathbb{C}\llangle x,y\rrangle}
\newcommand{\coh}{\mathop{\rm coh}\nolimits}

\newcommand{\uCM}{\mathop{\underline{\rm CM}}\nolimits}
\newcommand{\End}{\mathop{\rm End}\nolimits}
\newcommand{\Hom}{\mathop{\rm Hom}\nolimits}
\renewcommand*\partial{\textnormal{\reflectbox{6}}}

\newcommand{\per}{\mathop{\rm{Perf}}\nolimits}
\newcommand{\Db}{\mathop{\rm{D}^b}\nolimits}
\renewcommand{\mod}{\mathop{\textnormal{mod}}\nolimits}

\newcommand{\rmA}{\mathrm{A}}
\newcommand{\Acon}{\rmA_{\textnormal{con}}}
\newcommand{\AconTwoZ}{\rmA_{\textnormal{con}}^{2\mathbb{Z}}}

\newcommand{\Aconder}{\rmA_{\textnormal{con}}^{\textnormal{der}}}
\newcommand{\Bconder}{\mathrm{B}_{\textnormal{con}}^{\textnormal{der}}}

\newcommand{\Bcon}{\mathrm{B}_{\textnormal{con}}}
\newcommand{\Ccon}{\mathrm{C}_{\textnormal{con}}}

\begin{document}
\title{The taming of $3$-fold flops}
\author{Michael Wemyss}
\address{Michael Wemyss, School of Mathematics and Statistics,
University of Glasgow, University Place, Glasgow, G12 8QQ, UK.}
\email{michael.wemyss@glasgow.ac.uk}
\begin{abstract}
This article overviews recent progress in understanding and classifying one of the most elementary of birational surgeries within algebraic geometry, namely $3$-fold flops, through the development of a noncommutative version of singularity theory.

\end{abstract}

\maketitle

\section{Overview.}
Carrots generally work better than sticks, at least where birational geometry is concerned. In that context brute force attacks are often (but not always) just too tricky, especially as the dimension increases, and so time spent teasing out the true essence or structure of a problem is usually better than say writing out ever increasing systems of complicated polynomial equations.  

The same observation is true, but perhaps less so, for finite dimensional algebras. Here, whilst carrots remain prized, simply put, sticks are now a much more viable option. The finiteness of the underlying vector space allows all manner of techniques to be applied, and from these techniques, simplicity and elegance can follow.

\medskip
With the above in mind, the point of this survey is to explain the recent  classification of smooth $3$-fold flops within birational geometry.  In this story, the carrot is the \emph{contraction algebra}, discovered in \cite{DW1,DW3}, which through a sequence of steps explained in \S\ref{subsec:ContAlg} and \S\ref{subsec:DWtrue} below reduces the classification of flops to an isomorphism problem in noncommutative finite dimensional algebras.  It is not immediately clear that this helps, given that isomorphism problems in finite dimensional algebras are also in general painful, however here the stick enters, albeit briefly. A \emph{very small} number of noncommutative analytic changes in co-ordinates, which in many ways simply mimic known moves within classical singularity theory \`{a} la  Arnold \cite{Arnold},  simplify the situation dramatically. This leaves behind, in \S\ref{subsec:NCnormalIntro}, only a small set of possible contraction algebras.  

The classification of flops is then there for the taking, but is contingent on one final geometric construction phase, explained in \S\ref{subsec:ConstructFlops}, where conjectural noncommutative contraction criteria makes a bold prediction on how the large families of flops needed to complete the classification \emph{must} be constructed.  It cannot be said that the literature is overflowing with examples of flops, so thankfully this prediction turns out to be true, and in one stroke it uncovers large new families of flops, gives the complete classification, and also unifies all previous constructions. The end result is \ref{thm:GWThm} below, where an elegant list of noncommutative normal forms, akin to the lists produced by Arnold \cite{Arnold}, Du Val \cite{DuVal1,DuVal2} and Klein \cite{Klein}, provides the full and final classification.

\subsection{Language and setting.}
Both classification and in general the existence of contractions are necessarily local, and so throughout we consider a $3$-fold flopping contraction $\scrX\to\Spec\scrR$, for suitably local $\scrR$, and its flop $\scrY\to\Spec\scrR$, all of which is recalled in \S\ref{subsec:2Flops} (see also \cite{Kollar}). For single-curve flops, namely when there is only a single $\mathbb{P}^1$ above the origin, the cartoon picture is the following.
\begin{equation}
\begin{array}{c}
\begin{tikzpicture}
\node at (-1.25,0) {$\begin{tikzpicture}[scale=0.5]
\coordinate (T) at (-2.1,2);
\coordinate (B) at (-1.9,1);
\draw[red,line width=1pt] (T) to [bend left=25] (B);
\draw[color=blue!60!black,rounded corners=7pt,line width=1pt] (-0.5,0.1) -- (-1.75,0.3)-- (-3.4,0.2) -- (-3.8,1.5)-- (-3.6,2.8) -- (-1.75,3) -- (-0.75,2.8) -- (-0.3,2)-- cycle;
\end{tikzpicture}$};
\node at (1.25,0) {$\begin{tikzpicture}[scale=0.5]
\coordinate (T) at (2.1,2);
\coordinate (B) at (1.9,1);
\draw[red,line width=1pt] (T) to [bend right=25] (B);
\draw[color=blue!60!black,rounded corners=7pt,line width=1pt] (0.5,0.1) -- (1.75,0.3)-- (3.4,0.2) -- (3.8,1.5)-- (3.6,2.8) -- (1.75,3) -- (0.75,2.8) -- (0.3,2)-- cycle;
\end{tikzpicture}$};
\node at (0,-2) {\begin{tikzpicture}[scale=0.5]
\filldraw [red] (2.15,0.75) circle (1pt);
\draw[color=blue!60!black,rounded corners=5pt,line width=1pt] 
(0.5,0) -- (2,0.15)-- (4,0) -- (4.3,0.75)
-- (4,1.6) -- (2,1.35) -- (0.5,1.6) -- (0.2,0.75)-- cycle;
\end{tikzpicture}};
%%%
\draw[->, color=blue!60!black] (-0.5,-0.85) -- (-0.25,-1.5);
\draw[->, color=blue!60!black] (0.5,-0.85) -- (0.25,-1.5);
\node at (-2.5,0) {$\scrX$};
\node at (2.5,0) {$\scrY$};
\node at (-2,-2) {$\Spec\scrR$};
\end{tikzpicture}
\end{array}\label{eqn:SinglCurveFlop}
\end{equation}
Pausing the discussion of the presence or otherwise of a technical class of singularities on $\scrX$, for simplicity we will now assume that $\scrX$ (equivalently, $\scrY$) is smooth, and return to the presence of singularities in  \S\ref{subsec:Singular}.

Roughly speaking, flops can be viewed in two different ways. The first is that flops are the inevitable consequence of contracting curves, and as such they should be viewed as a $3$-fold `lifting' of various results of Artin~\cite{Artin} where chains of ($-2$)-curves are contracted.   The second is that flops are the answer to the lack of uniqueness in the output of the minimal model program, given that they are precisely the surgery that connects minimal models~\cite{KawaConnects}. The first viewpoint of flops is local, and leads to both classification and questions around contractibility, whereas the second is more global, and leads to the many connections with other parts of mathematics. Both viewpoints are important. This survey focuses almost entirely on the local, with the applications being relegated to a brief discussion within \S\ref{subsec:Global}, and references therein.

\subsection{What data to classify?}\label{subsec:WhatToClassify}
The single-curve cartoon picture in \eqref{eqn:SinglCurveFlop} is simultaneously illuminating and misleading. The picture is drawn to be deliberately symmetric, so that it is impossible to tell the two sides apart.  As such, if we were to `classify flops', it is going to be difficult to aim specifically for $\scrX\to\Spec\scrR$ without accidentally constructing $\scrY\to\Spec\scrR$.  Making the choice to aim for $\scrX$ also seems rather artificial, when $\scrY$ is just as good. The solution is to instead simply classify the possible $\scrR$ that can appear, as the two sides of \eqref{eqn:SinglCurveFlop} can then be obtained via blowing up the two generators $\pm 1$ of the class group $\mathrm{Cl}(\scrR)\cong\mathbb{Z}$, a canonical operation.

At this stage, considering the multi-curve situation becomes helpful. A given $2$-curve flopping contraction $\scrX\to\Spec\scrR$ still has a flop $\scrY\to\Spec\scrR$, under which both curves are flopped together.  As before, neither $\scrX$ or $\scrY$ can be preferred, and so aiming for $\scrR$ again seems sensible.  However, the situation is now complicated by the fact that above $\Spec\scrR$ sit additional flopping contractions. Indeed, labelling the curves in $\scrX$ by $\Curve_1$ and $\Curve_2$, it is possible to flop only curve $\Curve_1$ to obtain a new variety $\scrX_1$, or to flop only curve $\Curve_2$ to obtain a new variety $\scrX_2$. The outcome of these operations is summarised in the following picture.
\[
\begin{array}{c}
\begin{tikzpicture}
%%%LEFT RESOLUTION
\node at (-2.5,0) {$\begin{tikzpicture}[scale=0.5]
\coordinate (T1) at (2.1-4,2.5);
\coordinate (T) at (1.85-4,1.75);
\coordinate (TM) at (1.92-4,1.6);
\coordinate (B) at (2.1-4,0.8);
\node at (-1.6,2.1) {$\scriptstyle \Curve_1^+$};
\node at (-1.5,1.2) {$\scriptstyle \Curve_2$};
\draw[red,line width=1pt] (T1) to [bend right=25] (TM);
\draw[line width=1pt,red] (T) to [bend left=25] (B);

\draw[color=blue!60!black,rounded corners=7pt,line width=1pt] (-0.3,0.3) -- (-1.75,0.3)-- (-3.2,0.1) -- (-3.8,1.5)-- (-3.7,2.8) -- (-1.8,3.2) -- (-0.75,2.7) -- (-0.3,2)-- cycle;
\end{tikzpicture}$};
%%%MIDDLE RESOLUTION
\node at (0,0) {$\begin{tikzpicture}[scale=0.5]
\coordinate (T1) at (1.7-4,2.5);
\coordinate (T) at (1.9-4,1.7);
\coordinate (TM) at (1.92-4,1.6);
\coordinate (B) at (2.1-4,0.8);
\node at (-1.7,2.1) {$\scriptstyle \Curve_1$};
\node at (-1.5,1.2) {$\scriptstyle \Curve_2$};
\draw[red,line width=1pt] (T1) to [bend left=25] (TM);
\draw[line width=1pt,red] (T) to [bend left=25] (B);

\draw[color=blue!60!black,rounded corners=7pt,line width=1pt] (-0.5,0.1) -- (-1.75,0.3)-- (-3.4,0.2) -- (-3.8,1.5)-- (-3.6,2.8) -- (-1.75,3) -- (-0.75,2.8) -- (-0.3,2)-- cycle;
\end{tikzpicture}$};
%%%RIGHT
\node at (2.5,0) {$\begin{tikzpicture}[scale=0.5]
\coordinate (T1) at (1.7-4,2.5);
\coordinate (T) at (1.975-4,1.75);
\coordinate (TM) at (1.92-4,1.6);
\coordinate (B) at (1.8-4,0.8);
\node at (-1.7,2.1) {$\scriptstyle \Curve_1$};
\node at (-1.8,1.2) {$\scriptstyle \Curve_2^+$};
\draw[red,line width=1pt] (T1) to [bend left=25] (TM);
\draw[line width=1pt,red] (T) to [bend right=25] (B);

\draw[color=blue!60!black,rounded corners=7pt,line width=1pt] (-0.5,0.1) -- (-1.75,0.1)-- (-3.2,0.4) -- (-3.8,1.5)-- (-3.6,2.8) -- (-1.75,3) -- (-0.45,2.8) -- (-0.1,2)-- cycle;
\end{tikzpicture}$};
%%%
\node at (0,-2) {\begin{tikzpicture}[scale=0.6]
\filldraw [red] (2.15,0.75) circle (1pt);
\draw[color=blue!60!black,rounded corners=5pt,line width=1pt] 
(0.5,0) -- (2,0.15)-- (4,0) -- (4.3,0.75)
-- (4,1.6) -- (2,1.35) -- (0.5,1.6) -- (0.2,0.75)-- cycle;
\end{tikzpicture}};
%%%
\draw[->, color=blue!60!black] (0,-0.8) -- (0,-1.5);
\draw[->, color=blue!60!black] (-1.75,-0.8) -- (-1.25,-1.5);
\draw[->, color=blue!60!black] (1.75,-0.8) -- (1.25,-1.5);
%%%
\node at (0,1) {$\scrX=\scrX_0$};
\node at (2.5,1) {$\scrX_2$};
\node at (-2.5,1) {$\scrX_1$};
\node at (-2,-2) {$\Spec\scrR$};
\end{tikzpicture}
\end{array}
\]
Now the picture is not symmetric, and in general none of $\scrX_0$, $\scrX_1$ or $\scrX_2$ are isomorphic. Regardless, the process continues, iteratively flopping the new curves, and through this we obtain many smooth flopping contractions $\scrX_i\to\Spec\scrR$, $i=0,1,2,\hdots$. To illustrate the complexity, even understanding simply the number of contractions obtained via this iterative process in the $2$-curve case is complicated: by \cite{DW3, IW9} the answer is either $6, 8, 10, 12$, or $16$, numbers that are not immediately guessable.

\medskip
Given the above, we return to the general case of an $n$-curve flopping contraction $\scrX\to\Spec\scrR$, and ask the question of what to classify.  As for single- or two-curve flops, we cannot aim directly for $\scrX\to\Spec\scrR$, and aiming only for $\scrR$ seems to create uniqueness issues since it is no longer obvious how to easily distinguish or construct the original $\scrX\to\Spec\scrR$. The key point is that any specific flopping contraction \emph{cannot} be distinguished amongst the others, and so we should stop trying.  It turns out that all the $\scrX_i$ can still be constructed via blowing up regions inside the class group $\mathrm{Cl}(\scrR)\cong\mathbb{Z}^{\oplus n}$ \cite{Pinkham}, and as such, constructing all the $\scrX_i$ is both just as easy as, and also perhaps more canonical than, constructing only one. 

As a consequence, the initial aim of our classification is to specify all the possible $\scrR$ that can appear as the base of smooth flopping contractions, with the understanding that from that information, all the smooth flopping contractions above $\Spec\scrR$ can be obtained canonically. The theoretical aspects, particularly \S\ref{subsec:DWtrue} below, proceed on that basis.  However, from a more practical perspective later there is one significant downside to this approach, namely that even for single-curve flops the equations of the base $\scrR$ are not currently human-readable (see e.g.\ \cite[\S5.2]{BW4}), and it is not clear whether they ever will be.  This motivates the need for a more elegant form in which to encode and store the classification data, one from which $\scrR$ can be recovered in some non-explicit manner, e.g.\ through appealing to Mather--Yau \cite{MatherYau}, or via some abstract contraction theorem.   

The thirst for this more readable data is quenched by the (super)potential of the contraction algebra, explained in \S\ref{subsec:NCnormalIntro} below, which takes the form of an element in the free algebra in two variables.  This form has the benefit of being particularly elegant, and of being curiously similar to commutative singularity theory, but it also has a downside, namely of being almost the optimally worst form for geometers wishing to test geometric statements.  In part, this motivates the description of flops via gluing in \S\ref{subsec:ConstructFlops}, since this alternative description gives classification via clear geometric data.  It is our intention, through time, to give the classification in as many ways as possible, and to produce a $3$-fold version of the fifteen characterisations \`{a} la  \cite{Durfee}.

\subsection{The contraction algebra.}\label{subsec:ContAlg}
The \emph{contraction  algebra} $\Acon$ can be constructed in various ways, either as a noncommutative deformation algebra via deformation theory (see \cite{DW1,DW3}), or as a stable endomorphism ring inspired by cluster theory, or (in the smooth case here) as the Jacobi algebra of a quiver with potential. 

Below we will focus mostly on the cluster construction, since that is heavily used in all approaches to theoretical classification, then on the quiver with potential construction, since that is heavily used in the practical classification. We do however briefly pause to explain the importance of the deformation theoretic aspects, since it is through these that we are able to recover all known curve invariants from $\Acon$.  Indeed one example of this, of which there are many others, is Toda's dimension formula \cite{TodaGV}, which to simplify briefly and restrict to single-curve flops states that
\[
\dim_{\mathbb{C}}\Acon = \underbrace{\dim_{\mathbb{C}} \Acon^{\rm ab}}_{n_1} + \sum_{i=2}^{\ell} n_i\cdot i^2,
\]
where $\ell$ is Koll\'ar's \emph{length}, which is determined by Reid's general elephant \cite{Pagoda}, and the $n_i> 0$ are the Gopakumar--Vafa invariants \cite{KatzGV}.  Other numerical invariants are extracted using various other methods.  The point is, however, that $\Acon$ is strictly finer than any numerical curve invariant. The following was an early result.
\begin{theorem}[\cite{BW}]\label{thm:BW_GV_fail}
There are two single-curve flopping contractions $f_i\colon \scrX_i\to\Spec \scrR_i$, $i=1,2$, for which the following statements hold.
\begin{enumerate}
\item $\scrR_1$ is not analytically (or algebraically) isomorphic to $\scrR_2$.
\item The Gopakumar--Vafa invariants associated to $f_1$ and $f_2$ are the same.
\item The contraction algebras associated to $f_1$ and $f_2$ are not isomorphic.
\end{enumerate}
\end{theorem}
The contraction algebra is also finer than various motivic versions of curve-counting invariants \cite{Okke}.  It is thus not unreasonable that, as the finest curve invariant, $\Acon$ could be the classifying object.  The point to emphasise is that it is the algebra structure that gives the fine control: there is only one number $9$, and thus only one $9$-dimensional vector space, but there is more than one $9$-dimensional finite dimensional algebra.  From this, \ref{thm:BW_GV_fail} more-or-less immediately follows.

\medskip
We return to the construction of $\Acon$ as a cluster-tilted algebra, since this by design will also account for the lack of uniqueness of the $\scrX_i\to\Spec\scrR$ in \S\ref{subsec:WhatToClassify} above.  Given any fixed flopping contraction $\scrX\to\Spec\scrR$, following work of Bridgeland \cite{BridgelandFlops}, Van den Bergh \cite{VdB1d} constructs a tilting bundle $\scrO\oplus\scrN$ on $\scrX$, which after setting $N\colonequals \mathrm{H}^0(\scrN)$ induces a derived equivalence
\begin{equation}
\mathbf{R}\mathrm{Hom}_\scrX(\scrO\oplus\scrN,-)\colon\Db(\coh \scrX)\xrightarrow{\sim}\Db(\mod\End_{\scrR}(\scrR\oplus N)).\label{VdBtilt}
\end{equation}
From here, the contraction algebra is defined to be $\Acon\colonequals \underline{\End}_\scrR(N)$, the quotient ring of $\End_\scrR(N)$ by the two-sided ideal of all morphisms $N\to N$ that factor through a free $\scrR$-module.  This definition involves choice, since there are many different tilting bundles on $\scrX$ and thus equivalences \eqref{VdBtilt}, but the deformation theory construction and viewpoint (not explained here) shows that $\Acon$ depends only on $\scrX$.

The benefit of the above definition of $\Acon$ is that, via standard noncommutative resolution arguments, the module $N$ is a cluster-tilting object when viewed as an object inside the singularity category of $\scrR$, and cluster theory is both well developed and good at iterating.  In particular, it is not unreasonable to ask what the other cluster-tilting objects are, and what they encode.

\begin{theorem}[\cite{HomMMP}]\label{thm:AMCorr}
For any $\scrR$ which is the base of a $3$-fold smooth flopping contraction, there is a bijection between the set of smooth flopping contractions $\scrX\to\Spec\scrR$ and the set of cluster-tilting objects in the singularity category of $\scrR$.
\end{theorem}

To each smooth flopping contraction $f_i\colon \scrX_i\to\Spec\scrR$, the bijection associates a cluster-tilting object $N_i$, in such a way that $\underline{\End}_\scrR(N_i)$ is the contraction algebra for $f_i$. Now, if the contraction algebra is to be the classifying object, following the logic in \S\ref{subsec:WhatToClassify} we cannot aim for any specific $\underline{\End}_\scrR(N_i)$, since that would make $\scrX_i$ distinguished.  The various $\underline{\End}_\scrR(N_i)$ are not isomorphic in general, and here in the finite dimensional algebra world this lack of uniqueness is more acute, since the previous solution of aiming for a common `base'  $\scrR$ cannot easily be replicated.  This is solved by the following result of August, which asserts that the contraction algebras $\underline{\End}_\scrR(N_i)$ are all derived equivalent and thus cannot be homologically distinguished, but also much more remarkably, that there are no other algebras in their derived equivalence class. 
\begin{theorem}[\cite{August1,August2}]\label{thm:Jenny}
Given a smooth $3$-fold flopping contraction $\scrX\to\Spec\scrR$ with associated $\Acon$, and an arbitrary basic algebra $\mathrm{B}$, the following conditions are equivalent.
\begin{enumerate}
\item $\Db(\mod\Acon)\simeq \Db(\mod \mathrm{B}) $.
\item $\mathrm{B}\cong\Ccon$ for a contraction algebra $\Ccon=\underline{\End}_\scrR(N_i)$ associated to one of the (finitely many) smooth flopping contractions $\scrX_i\to\Spec \scrR$. 
\end{enumerate}
\end{theorem}
The upshot, as in \S\ref{subsec:WhatToClassify}, is that we thus stop trying to aim for a specific $\underline{\End}_\scrR(N_i)$ and instead aim for the whole derived equivalence class. 

\subsection{Theoretical classification.}\label{subsec:DWtrue} With the above said, classification becomes the only thing that it can. Namely, that the contraction algebra classifies, and it does so in such a way that the (inevitable) non-uniqueness issues are accounted for in the most direct manner.

\begin{theorem}[previously the Donovan--Wemyss Conjecture]\label{conj:DW}
Suppose that $\scrX_1\to\Spec \scrR_1$ and $\scrX_2\to\Spec \scrR_2$ are two smooth $3$-fold flopping contractions, where both $\scrR_i$ are complete local.  Denote their corresponding contraction algebras by $\Acon$ and $\Bcon$ respectively.  Then $\scrR_1\cong \scrR_2$ if and only if $\Db(\mod\Acon)\simeq\Db(\mod\Bcon)$.
\end{theorem}

Using \ref{thm:Jenny} it is possible to replace the derived equivalence in \ref{conj:DW} with some isomorphism problem, but it should be remarked that this is not in general as simple as asking for an isomorphism between the two given $\Acon$ and $\Bcon$.  As such, this replacement comes at the cost of introducing a third $\Ccon$, which then makes the statement both less symmetric, and less appealing.  In the single-curve setting, which is a key benchmarking case below (and also the most important, since any flop decomposes as a sequence of single-curve flops), contraction algebras are local and so there \ref{conj:DW} does reduce to the more elementary statement that $\scrR_1\cong \scrR_2$ if and only if $\Acon\cong\Bcon$. 

Regardless, the main point of \ref{conj:DW} is the reduction in complexity, as it asserts that $3$-fold geometry is classified by some finite dimensional data, at the cost of going noncommutative. 

\medskip
We now briefly sketch the proof of \ref{conj:DW}, leaving most details to \S\ref{subsec:Derived}. The paper \cite{DW1} plus \ref{thm:Jenny} proves ($\Rightarrow$), and so it suffices to prove ($\Leftarrow$).  The overall strategy relies on an observation of Kalck--Yang \cite{KalckYang}, and of Buchweitz at MSRI in 2013.  The contraction algebra $\Acon$ can be viewed as the cohomology $\mathrm{H}^0(B)$, where $B=\Aconder$ is a DG-algebra now known as the derived contraction algebra (see also Booth \cite{Booth}).  It is tautologically true that $B$ determines both the category $\per(B)$ and the category $\mathrm{D}_{\mathrm{fd}}(B)$.  The fact that $\scrX$ is smooth translates into the existence of a noncommutative crepant resolution (NCCR), and from this Kalck--Yang \cite{KalckYang} show that there is a natural inclusion $\mathrm{D}_{\mathrm{fd}}(B)\subset \per(B)$, and further that the quotient category is the classical singularity category $\uCM\scrR$.  The upshot is that $B$ determines the triangulated category $\uCM\scrR$, and since $\scrR$ is a hypersurface, this is very close to determining $\scrR$.

The first part of the following result asserts that with some extra work this triangulated statement can be upgraded to the DG-enhancement, in which case by Hua--Keller \cite{HuaKeller} we have access to the classical Tjurina algebra (via Hochschild cohomology), and thus the Mather--Yau theorem \cite{MatherYau}.

\begin{theorem}\label{thm:DWThm}
For a smooth $3$-fold flopping contraction $\scrX\to\Spec\scrR$, the following statements hold.
\begin{enumerate}
\item\label{thm:DWThm 1} \textnormal{\cite{Booth, HuaKeller}} The derived contraction algebra $\Aconder$ determines $\scrR$.
\item\label{thm:DWThm 2}  \textnormal{\cite{JKM,KLW}} $\Acon$ determines the derived contraction algebra $\Aconder$, and thus \ref{conj:DW} is true.
\end{enumerate}
\end{theorem}

The second part of \ref{thm:DWThm} is the most involved of all, and we sketch this in \S\ref{subsec:Derived} later.

\subsection{Towards noncommutative normal forms.}\label{subsec:NCnormalIntro}
Even in the single-curve case, the uncharitable reading of \ref{thm:DWThm} is that those finite dimensional algebras $\Acon$ which come from geometry classify geometry. If it is not possible to recognise which finite dimensional algebras are isomorphic to some $\Acon$, has the process of passing from $3$-fold geometry to finite dimensional algebras simply shifted the complexity elsewhere?

At this stage, it is instructive to recall another incarnation of the contraction algebra, namely as a Jacobi algebra $\Jac(g)$ of a quiver with potential \cite{VdBCY}.  We now restrict the discussion to the single-curve case, to simplify and also since for classification purposes this is the most important case.  Consider the free algebra in two variables $\mathbb{C}\langle x,y\rangle$, and its complete local version $\ringtwo$, where infinite sums are permitted. As explained in e.g.\ \cite{DWZ}, for an element $g\in\ringtwo$, the Jacobi algebra is defined to be
\[
\Jac(g) = \frac{\ringtwo}{\lcl\updelta_{x}g ,\updelta_{y}g\rcl}
\]
where $\updelta_{x}$ and $\updelta_y$ are the cyclic derivatives, and $\lcl -\rcl$ denotes the closure of the two-sided ideal.

One property contraction algebras of flops satisfy is that they are finite dimensional \cite{DW1}. Another property that they enjoy, essentially due to Reid's general elephant \cite{Pagoda}, is that they admit a central element $z\in\scrZ(\Jac(g))$ such that $\Jac(g)/\lcl z\rcl$ is isomorphic to an algebra of the form $e\Uppi e$, where $\Uppi$ is a preprojective algebra of ADE type, and $e$ is a primitive idempotent (see e.g.\ \cite{DW1}).  This property is very strong, since given an arbitrary Jacobi algebra, it is not even clear whether the centre is non-trivial.

Given this, and referring the reader to \ref{def:alglengthintro} below for more details, we say that an arbitrary algebra $B$ has \emph{algebraic length} $\ell$ if it admits a central element as above, where the idempotent $e$ has associated number $\ell$ in \eqref{eqn:allDynkin}.  The following result, which is rather surprising, shows that the algebraic length is well-defined, and puts the options into six families.  Both the statement and proof remain entirely within algebra, and both rely only on a series of rather elementary noncommutative analytic changes in coordinates.
\begin{theorem}[\cite{BW5}]\label{thm:BW_alg_length}
Let $g\in\ringtwo$. If $\dim_\mathbb{C}\Jac(g)<\infty$, and further $g$ has an algebraic length, then $\Jac(g)$ is isomorphic to $\Jac(h)$ for some $h$ in the following list
\[
\begin{array}{l}
\begin{tabular}{p{0.3cm}p{3.5cm}p{4cm}}
$1$ & $x^2 + y^n$ & $n\geq 2$, \\[2pt]
$2$ &$xy^2+ x^{2n}+\upvarepsilon x^{2m-1}$ & $m,n\geq2$, $\upvarepsilon\in\{0,1\}$ \\[2pt]
$3$ &$x^3 + xy^3 + p(y)$ & $0\neq p\in\mathbb{C}[y],\, \mathrm{ord}(p)\geq4$ \\[2pt]
\end{tabular}\\
\begin{tabular}{p{0.3cm}p{5cm}p{1cm}}
$4$ & $x^3 + y^4 + xyxy^2 + \scrO_6$\\[2pt]
$5$ & $x^3 + xyxy + y^5 \,\,\,+ \scrO_6$\\[2pt]
$6$ & $x^3 + y^4 + (xy)^3 + \uplambda(xy^2)^2 +\, \scrO_7$&$\uplambda\neq 0$
%\bottomrule
\end{tabular}
\end{array}
\]
where the algebraic length $\ell$ indexes the rows, and $\scrO_n$ indicates possible higher order terms of degree $n$ and higher. 
\end{theorem}
Since contraction algebras of smooth single-curve flops satisfy the assumptions of \ref{thm:BW_alg_length} they must, at worst, be a subset of the Jacobi algebras listed.  The classification of smooth $3$-fold flops would thus follow at once, if all the $\Jac(h)$ for $h$ in the above list are in fact contraction algebras. 

\subsection{Construction of flops.}\label{subsec:ConstructFlops}
The list in \ref{thm:BW_alg_length} helps, to the extent that it makes a prediction.  The mild downside is that whole new families of flops need to be constructed, both for the fancier reason that we wish to demonstrate that everything in \ref{thm:BW_alg_length} is isomorphic to a contraction algebra, but also from the much easier reason that if a family of flops exists, it necessarily needs to be in the classification.

Finding or constructing flops is hard, partly since it has a needle-in-a-haystack feel to it; the world of all curves is wild, and the subset of curves that flop are only a tiny part of it.  This problem can be solved, if we allow ourselves to be motivated by a different noncommutative conjecture, namely the contraction conjecture of \cite{BW3, Kawamata}. This asserts that curves \emph{will} contract (suitably locally), provided that their noncommutative deformation algebra is finite dimensional. What follows gives serious evidence for this conjecture, however what is perhaps the most remarkable aspect of the whole story is just how practically useful the noncommutative structure is.

\medskip
For any choice of finitely many scalars $\uplambda_{jk}$ where $j,k\geq 0$ and $j+k\geq 2$, as in \cite{Ferrari, KatzVersal} consider the variety $\scrX$ defined by the two open  patches $\scrU_1=\mathbb{A}^3_{a,v_2,v_1}$ and $\scrU_2=\mathbb{A}^3_{b,w_2,w_1}$ glued via the rule
\begin{equation}\label{eq!glue!intro}
(a,v_2,v_1)\stackrel{a\neq 0}{\longleftrightarrow}(a^{-1}, \,a^{-1}v_2, \,\,a^{3}v_1+\partial E\,), \quad\textnormal{where}\quad \partial E \colonequals \sum_{j+k\geq 2} \uplambda_{jk}a^{2-k}v_2^{j+k-1}.
\end{equation}
The locus $v_2=v_1=0$ glues with the locus $w_2=w_1=0$ to create a single $\mathbb{P}^1\subset\scrX$.  Given this setup, the paper \cite{BW2} computes the noncommutative deformation algebra of $\scrO_{\mathbb{P}^1}\in\coh\scrX$, which as it turns out is presented by some formula depending on the $\uplambda_{jk}$, the details of which are not important here.

The key observation is the following. Given any $h$ in the list in  \ref{thm:BW_alg_length}, which we emphasise is an element of $\ringtwo$ and thus noncommutative data, using the above construction it is easy to produce $\Curve=\mathbb{P}^1\subset\scrX$ in which the noncommutative deformation algebra is isomorphic to $\Jac(h)$.  If you believe the contraction conjecture, then it simply \emph{must} be the case that the curve $\Curve$ locally contracts to give a flopping contraction, since all the $\Jac(h)$ in \ref{thm:BW_alg_length} are finite dimensional.

This prediction is akin to locating the needles in the haystack: once their location has been predicted, verifying that they are needles is not the difficult part. Indeed, from here all noncommutative technology melts safely into the background, and using only techniques from algebraic geometry that are over 30 years old \cite{Jimenez}, the following asserts that indeed the predicted curves contract, and thus build all flops. 
\begin{theorem}[\cite{BW6}]\label{thm:GeomReal}
Consider the variety $\scrX$ obtained by gluing in \eqref{eq!glue!intro}, in which $\partial E$ is one of the following.
\begin{center}
\begin{tabular}{cllp{4cm}}
$1$ &$ a^2z^{\phantom 2} + a^{2-n}z^{n-1}$& $n\geq 2$\\[3pt]

$2$ &${\phantom a^{\phantom 2}}z^2 + a^2z^{2n-1} + \upvarepsilon a^2z^{2m-2} $ &  
$m,n\geq2$, $\upvarepsilon\in\{0,1\}$\\[3pt]

$3$&$a^2z^2 + a^{-1}z^3 + \sum_{k\geq 3}\upalpha_k a^{2-k}z^k$ &\mbox{some }$\upalpha_j\neq 0$ \\[3pt]

$4$&$a^2z^2 + a^{-2}z^3 + a^{-1}z^4 + \scrO_5 $ & \\  [3pt]

$5$&$a^2z^2 + \phantom{a^{-2}}z^3 + a^{-3}z^4 + \scrO_5$ &\\[3pt]

$6$ & $a^2z^2 +  a^{-2}z^3 + a^{-1}z^5 + \uplambda a^{-2}z^5 + \scrO_6$  & $\uplambda\in\mathbb{C}^*$ 
\end{tabular}
\end{center}
Then $\Curve=\mathbb{P}^1\subset\scrX$ locally contracts, creating a smooth flopping contraction $\hat{\scrX}\to\Spec \scrR$ in which the contraction algebra realises the Jacobi algebra predicted by \ref{thm:BW_alg_length}.
\end{theorem}
Before \ref{thm:GeomReal} there were only four examples of higher length flops ($\ell\geq 3$) known in the literature \cite{Karmazyn}. This gave the impression that $3$-fold flops should behave like ADE surfaces, in the sense that the big families are Type A and D ($\ell=1,2$ here) with only a few sporadic cases in Type E ($\ell\geq 3$). The above \ref{thm:GeomReal} in fact shows that the opposite is true: for $3$-fold flops the small countable cases are actually $\ell=1,2$, whereas the $\ell\geq 3$ cases live in families, all of which are uncountable \cite{BW5}.

\subsection{The final classification.}
With the theoretical classification of \ref{thm:DWThm} in place, together with the algebraic restriction of possible normal forms of \ref{thm:BW_alg_length} and their geometric realisation in  \ref{thm:GeomReal}, all the ingredients are in place to state the final classification. 

To briefly recap the overall strategy: after reducing the problem to an isomorphism problem in Jacobi algebras in \ref{thm:DWThm}, the development of a noncommutative version of singularity theory \`{a} la  Arnold \cite{Arnold} gives the needed algebraic restriction of possible normal forms in \ref{thm:BW_alg_length}. This is the key reduction step. How to build the flops is then \emph{predicted} from this noncommutative data, via the contraction conjecture and specific gluings in \eqref{eq!glue!intro}. Verifying this prediction uses only geometric techniques, and is comparatively straightforward.  The following is the main result.

\begin{theorem}\label{thm:GWThm}
Smooth single-curve $3$-fold flops are classified by the following equivalent information.
\begin{enumerate}
\item The isomorphism classes of the Jacobi algebras $\Jac(h)$ in \ref{thm:BW_alg_length}.
\item The isomorphism classes of the gluing data $\partial E$  in \ref{thm:GeomReal}.
\end{enumerate}
\end{theorem}
With some extra work, again in the spirit of classical singularity theory, it is possible to improve the noncommutative forms within \ref{thm:BW_alg_length} so that there are no isomorphisms between different elements.  After this step, the `isomorphism classes' part of  \ref{thm:GWThm} can be removed.  Of course, this step is complicated by the existence of moduli, and details will appear elsewhere.

\subsection{Other related problems.} \label{subsec:Global}
This article outlines the straight-line approach to classification and deliberately omits many of the wider structures, used both technically and conceptually, that also turn out to facilitate applications to other areas.   These include Coxeter theory, commutative algebra and modifying modules, derived categories, stability conditions, faithful actions, symplectic topology, the $K(\pi,1)$-conjecture, Lie theory, and many others, where we refer the reader to e.g.\ \cite{lockdown} for more details.  The noncommutative singularity aspects are explained in more depth in the survey \cite{BW4}. 

\subsection{Singular classification.}\label{subsec:Singular}
Singularities are a fact of life in the minimal model program, and much of what is said within this article should also in time extend to the case of singular flops $\scrX\to\Spec\scrR$ where $\scrX$ has Gorenstein terminal singularities. Certainly all the applications listed within \S\ref{subsec:Global} now all work within this broader framework, but the main issue for classification is how to replace the noncommutative singularity theory aspects, since superpotentials are so heavily tied to the smooth CY setting.

It should be remarked that \ref{thm:DWThm}\eqref{thm:DWThm 1} is still true in the singular setting, but as in the smooth case, it is not expected that the full derived structure is necessary. The question remains as to what is the \emph{minimal} amount of information needed for classification when $\scrX$ is singular, and \cite{JKM, KLW} (which solve \ref{thm:DWThm}\eqref{thm:DWThm 2}) are both conspicuously quiet in this regard.

\section{Flops and contraction algebras.}

This section recalls some basic facts and definitions in \S\ref{subsec:2Flops}, discusses the key length invariant used for classification in \S\ref{subsec:length}, before sketching the proof of \ref{thm:DWThm}\eqref{thm:DWThm 2} in \S\ref{subsec:Derived}.

\subsection{Flops, length and slicing.}\label{subsec:2Flops}
There are various equivalent definitions of flops in the literature (see e.g.\ \cite{KMM, Kollar}). For our purposes, the following suffices.

\begin{definition}
Suppose that $f\colon \scrX\to \Spec \scrR$ is a crepant projective birational morphism, where $\scrR$ is isolated complete local \textnormal{cDV}.  Then $f^+\colon \scrX^+\to \Spec\scrR$ is called the flop of $f$ if for every line bundle $\scrL=\scrO_\scrX(D)$ on $\scrX$ such that $-D$ is $f$-nef, then the proper transform of $D$ is $\mathbb{Q}$-Cartier, and $f^+$-nef. 
\end{definition}
In the above,  a $\mathbb{Q}$-Cartier divisor $D$ is called $f$-nef if $D\cdot C\geq 0$ for all curves contracted by $f$.  The main point is that, under proper transform, the intersection of divisors with curves swaps sign.  Loosely speaking, we view the flop as a surgery in which the original curves in $\scrX$ above the origin $0$ are replaced by new curves in $\scrX^+$, in such a way that $\scrX\setminus f^{-1}(0)\xrightarrow{\sim} \scrX^+\setminus (f^+)^{-1}(0)$.

From the early days there has been a strong relationship between flops and ADE surface singularities. Indeed, Reid's general elephant \cite{Pagoda} asserts that for generic $g\in\scrR$, the base change diagram
\begin{equation}
\begin{array}{c}
\begin{tikzpicture}
\node (A) at (0,0) {$\scrY$};
\node (a) at (2,0) {$\scrX$};
\node (B) at (0,-1.25) {$\Spec\scrR/g$};
\node (b) at (2,-1.25) {$\Spec\scrR$};
\draw[->] (A)--(a);
\draw[->] (A)--(B);
\draw[->] (B)--(b);
\draw[->] (a)--(b);
\end{tikzpicture}
\end{array}\label{eqn:BaseChange}
\end{equation}
results in an ADE surface singularity $\scrR/g$, and $\scrY\to\Spec\scrR/g$ is a crepant partial resolution.  There are two main benefits: the first is that the above grants access to ADE-type combinatorics, and the second is that flops can be constructed, at least in principle, as a one-parameter deformation of ADE surface singularities.

The crepant partial resolution $\scrY\to\Spec \scrR/g$  in \eqref{eqn:BaseChange} is dominated by the minimal resolution $\scrZ\to\scrY\to\Spec\scrR/g$.  In particular, $\scrY$ can be obtained from $\scrZ$ by blowing down curves, and by McKay correspondence the curves in $\scrZ$ are in bijection with the nodes of an ADE Dynkin diagram. When $\scrX\to\Spec\scrR$ is a single-curve flop, $\scrY\to\Spec \scrR/g$ contains a single curve, thus $\scrY$ can be fully described by the pair $(\Delta,v)$ where $\Delta$ is the Dynkin diagram, and $v$ is the node corresponding to the unique curve which survives under the map $\scrZ\to\scrY$.  Katz--Morrison \cite{KM}, then later Kawamata \cite{KawaKM}, proved that for smooth single-curve flops only the following six pairs $(\Delta,v)$ can occur, where by convention we draw the $\Delta$, and $v$ is the vertex $\circ$. 
\[
\begin{tikzpicture}[scale=0.7]
\node at (0,0.1) {$
\begin{tikzpicture}[scale=0.45]
 \node (0) at (0,0) [DW] {};
  \node (1b) at (1,0.75) {};
 \node (0a) at (0,-0.5) {$\scriptstyle 1$};
 \node[draw=none] at (0,1) {};
\end{tikzpicture}$};
\node at (2,0) {\begin{tikzpicture}[scale=0.45]
 \node (0) at (0,0) [DB] {};
 \node (1) at (1,0) [DW] {};
 \node (1b) at (1,0.75) [DB] {};
 \node (2) at (2,0) [DB] {};
 \node (1a) at (1,-0.5) {$\scriptstyle 2$};
\draw [-] (0) -- (1);
\draw [-] (1) -- (2);
\draw [-] (1) -- (1b);
\end{tikzpicture}};
\node at (5.1,0) {$\begin{tikzpicture}[scale=0.45]
 \node (m1) at (-1,0) [DB] {};
 \node (0) at (0,0) [DB] {};
 \node (1) at (1,0) [DW] {};
 \node (1b) at (1,0.75) [DB] {};
 \node (2) at (2,0) [DB] {};
 \node (3) at (3,0) [DB] {};
 \node (1a) at (1,-0.5) {$\scriptstyle 3$};
\draw [-] (m1) -- (0);
\draw [-] (0) -- (1);
\draw [-] (1) -- (2);
\draw [-] (2) -- (3);
\draw [-] (1) -- (1b);
\end{tikzpicture}$};
\node at (9,-0) {$\begin{tikzpicture}[yscale=0.45,xscale=-0.45]
 \node (m1) at (-1,0) [DB] {};
 \node (0) at (0,0) [DB] {};
 \node (1) at (1,0) [DW] {};
 \node (1b) at (1,0.75) [DB] {};
 \node (2) at (2,0) [DB] {};
 \node (3) at (3,0) [DB] {};
 \node (4) at (4,0) [DB] {};
 \node (1a) at (1,-0.5) {$\scriptstyle 4$};
\draw [-] (m1) -- (0);
\draw [-] (0) -- (1);
\draw [-] (1) -- (2);
\draw [-] (2) -- (3);
\draw [-] (3) -- (4);
\draw [-] (1) -- (1b);
\end{tikzpicture}$};
\node at (3,-1.5) {$\begin{tikzpicture}[yscale=0.45,xscale=-0.45]
 \node (m1) at (-1,0) [DB] {};
 \node (0) at (0,0) [DB] {};
 \node (1) at (1,0) [DB] {};
 \node (1b) at (1,0.75) [DB] {};
 \node (2) at (2,0) [DW] {};
 \node (3) at (3,0) [DB] {};
 \node (4) at (4,0) [DB] {};
  \node (5) at (5,0) [DB] {};
\node (2a) at (2,-0.5) {$\scriptstyle 5$};
\draw [-] (m1) -- (0);
\draw [-] (0) -- (1);
\draw [-] (1) -- (2);
\draw [-] (2) -- (3);
\draw [-] (3) -- (4);
\draw [-] (4) -- (5);
\draw [-] (1) -- (1b);
\end{tikzpicture}$};
\node at (8,-1.5) {$\begin{tikzpicture}[yscale=0.45,xscale=-0.45]
 \node (m1) at (-1,0) [DB] {};
 \node (0) at (0,0) [DB] {};
 \node (1) at (1,0) [DW] {};
 \node (1b) at (1,0.75) [DB] {};
 \node (2) at (2,0) [DB] {};
 \node (3) at (3,0) [DB] {};
 \node (4) at (4,0) [DB] {};
  \node (5) at (5,0) [DB] {};
 \node (1a) at (1,-0.5) {$\scriptstyle 6$};
\draw [-] (m1) -- (0);
\draw [-] (0) -- (1);
\draw [-] (1) -- (2);
\draw [-] (2) -- (3);
\draw [-] (3) -- (4);
\draw [-] (4) -- (5);
\draw [-] (1) -- (1b);
\end{tikzpicture}$};
\end{tikzpicture}
\]
The above six cases can alternatively be restated using Koll\'ar's length invariant $\ell$.  Reading left to right, the cases are $\ell=1,\hdots,6$ respectively.

The upshot is that flops naturally split into six families according to the length invariant $\ell$, but really this is only data from surfaces.  Slicing from $3$-folds to surfaces results in a \emph{significant} loss of information, and aside from putting flops into convenient buckets, slicing does not allow for an obvious analytic $3$-fold classification.

\subsection{Algebraic length.}\label{subsec:length}
The outlook from the noncommutative perspective is much more optimistic. The point here is that it is possible to replicate  Koll\'ar's length in a purely algebraic form, one in which isomorphism classes of finite dimensional algebras can be readily manipulated.  

Leading into the next definition, recall the following ADE Dynkin diagrams, where the vertices have been labelled with the information of the highest root.
%Some text here
\begin{equation}\label{eqn:allDynkin}
\begin{array}{c}
\begin{tikzpicture}
\node at (0,0) {$A_n$};
%\node at (0,-0.2) {$\scriptstyle n\geq 1$};
\node at (2.3,0) {$
\begin{tikzpicture}[xscale=0.6]
 \node (0) at (0,0) [DW] {};
 \node (1) at (1,0) [DW] {};
 \node (2) at (2,0) [DW] {};
 \node (3) at (4,0) [DW] {};
 \node (4) at (5,0)[DW] {};
 \node at (3,0) {$\cdots$};
 \node (0a) at (0,-0.3) {$\scriptstyle 1$};
 \node (1a) at (1,-0.3) {$\scriptstyle 1$};
 \node (2a) at (2,-0.3) {$\scriptstyle 1$};
 \node (3a) at (4,-0.3) {$\scriptstyle 1$};
 \node (4a) at (5,-0.3) {$\scriptstyle 1$};
\draw [-] (0) -- (1);
\draw [-] (1) -- (2);
\draw [-] (2) -- (2.6,0);
\draw [-] (3.4,0) -- (3);
\draw [-] (3) -- (4);
\end{tikzpicture}$};
\node at (5,0) {$D_n$};
%\node at (5,-0.2) {$\scriptstyle n\geq 4$};
\node at (7.3,0) {\begin{tikzpicture}[xscale=0.6,yscale=0.7]
 \node (0) at (0,0) [DW] {};
 \node (1) at (1,0) [DW] {};
 \node (1b) at (1,0.75) [DW] {};
 \node (2) at (2,0) [DW] {};
 \node (3) at (4,0) [DW] {};
 \node (4) at (5,0)[DW] {};
 \node at (3,0) {$\cdots$};
 \node (0a) at (0,-0.3) {$\scriptstyle 1$};
 \node (1a) at (1,-0.3) {$\scriptstyle 2$};
 \node (1ba) at (1.3,0.75) {$\scriptstyle 1$};
 \node (2a) at (2,-0.3) {$\scriptstyle 2$};
 \node (3a) at (4,-0.3) {$\scriptstyle 2$};
 \node (4a) at (5,-0.3) {$\scriptstyle 1$};
\draw [-] (0) -- (1);
\draw [-] (1) -- (2);
\draw [-] (2) -- (2.6,0);
\draw [-] (3.4,0) -- (3);
\draw [-] (3) -- (4);
\draw [-] (1) -- (1b);
\end{tikzpicture}};
\node at (0,-1.5) {$\begin{tikzpicture}[yscale=0.7,xscale=0.6]
\node at (-0.5,0.75) {$E_6$};
 \node (m1) at (-1,0) [DW] {};
 \node (0) at (0,0) [DW] {};
 \node (1) at (1,0) [DW] {};
 \node (1b) at (1,0.75) [DW] {};
 \node (2) at (2,0) [DW] {};
 \node (3) at (3,0) [DW] {};
\node (m1a) at (-1,-0.3) {$\scriptstyle 1$};
 \node (0a) at (0,-0.3) {$\scriptstyle 2$};
 \node (1a) at (1,-0.3) {$\scriptstyle 3$};
 \node (1ba) at (1.3,0.75) {$\scriptstyle 2$};
 \node (2a) at (2,-0.3) {$\scriptstyle 2$};
 \node (3a) at (3,-0.3) {$\scriptstyle 1$};
\draw [-] (m1) -- (0);
\draw [-] (0) -- (1);
\draw [-] (1) -- (2);
\draw [-] (2) -- (3);
\draw [-] (1) -- (1b);
\end{tikzpicture}$};
\node at (3.5,-1.5) {$\begin{tikzpicture}[yscale=0.7,xscale=0.6]
\node at (-0.5,0.75) {$E_7$};
 \node (m1) at (-1,0) [DW] {};
 \node (0) at (0,0) [DW] {};
 \node (1) at (1,0) [DW] {};
 \node (1b) at (1,0.75) [DW] {};
 \node (2) at (2,0) [DW] {};
 \node (3) at (3,0) [DW] {};
 \node (4) at (4,0) [DW] {};
\node (m1a) at (-1,-0.3) {$\scriptstyle 2$};
 \node (0a) at (0,-0.3) {$\scriptstyle 3$};
 \node (1a) at (1,-0.3) {$\scriptstyle 4$};
 \node (1ba) at (1.3,0.75) {$\scriptstyle 2$};
 \node (2a) at (2,-0.3) {$\scriptstyle 3$};
 \node (3a) at (3,-0.3) {$\scriptstyle 2$};
 \node (4a) at (4,-0.3) {$\scriptstyle 1$};
\draw [-] (m1) -- (0);
\draw [-] (0) -- (1);
\draw [-] (1) -- (2);
\draw [-] (2) -- (3);
\draw [-] (3) -- (4);
\draw [-] (1) -- (1b);
\end{tikzpicture}$};
\node at (8.2,-1.5) {$\begin{tikzpicture}[yscale=0.7,xscale=0.6]
\node at (-0.5,0.75) {$E_8$};
 \node (m1) at (-1,0) [DW] {};
 \node (0) at (0,0) [DW] {};
 \node (1) at (1,0) [DW] {};
 \node (1b) at (1,0.75) [DW] {};
 \node (2) at (2,0) [DW] {};
 \node (3) at (3,0) [DW] {};
 \node (4) at (4,0) [DW] {};
  \node (5) at (5,0) [DW] {};
\node (m1a) at (-1,-0.3) {$\scriptstyle 2$};
 \node (0a) at (0,-0.3) {$\scriptstyle 4$};
 \node (1a) at (1,-0.3) {$\scriptstyle 6$};
 \node (1ba) at (1.3,0.75) {$\scriptstyle 3$};
 \node (2a) at (2,-0.3) {$\scriptstyle 5$};
 \node (3a) at (3,-0.3) {$\scriptstyle 4$};
 \node (4a) at (4,-0.3) {$\scriptstyle 3$};
 \node (5a) at (5,-0.3) {$\scriptstyle 2$};
\draw [-] (m1) -- (0);
\draw [-] (0) -- (1);
\draw [-] (1) -- (2);
\draw [-] (2) -- (3);
\draw [-] (3) -- (4);
\draw [-] (4) -- (5);
\draw [-] (1) -- (1b);
\end{tikzpicture}$};
\end{tikzpicture}
\end{array}
\end{equation}
To each of these graphs is an associated preprojective algebra $\Uppi$ (see e.g.\ \cite{BCB}),  which is finite dimensional.  The vertices correspond to primitive idempotents $e$ in $\Uppi$, thus each primitive idempotent of $\Uppi$ comes equipped with the number in \eqref{eqn:allDynkin} attached to its corresponding vertex.

\begin{definition}\label{def:alglengthintro}
An element $g\in\ringtwo$ is said to have \emph{algebraic length $\ell$}, where $1\leq \ell\leq 6$, if there exists a central element $z\in\Jac(g)$ such that $\Jac(g)/\lcl z\rcl \cong e\Uppi e$, for some \textnormal{ADE} preprojective algebra $\Uppi$ and some primitive idempotent $e$ which has associated number $\ell$. 
\end{definition}

It is certainly not clear that the algebraic length is well-defined. In general, there can exist central elements $z_1,z_2\in\Jac(g)$ such that $\Jac(g)/\lcl z_1\rcl \cong e_1\Uppi_1 e_1$ and $\Jac(g)/\lcl z_2\rcl \cong e_2\Uppi_2e_2$ for different \textnormal{ADE} preprojective algebras $\Uppi_1,\Uppi_2$ and different primitive idempotents $e_1$ and $e_2$.  Whilst this can happen, the point in the following result is that \emph{always} $\ell_1=\ell_2$, i.e.\ the numbers associated to the idempotents $e_1$ and $e_2$ in \eqref{eqn:allDynkin} are the same. This result is the single-most important step in the classification of flops, since it simultaneously achieves two things: it provides a purely algebraic version of Reid's general elephant, and also does so in such a way that the isomorphism classes of the Jacobi algebras (and thus, by \ref{thm:DWThm}\eqref{thm:DWThm 2}, the $3$-fold geometry) is both reduced and controlled.

\begin{theorem}[\cite{BW5}]
The algebraic length is well-defined. Furthermore, up to isomorphism of Jacobi algebras, all $g$ with algebraic length for which $\dim_\mathbb{C}\Jac(g)<\infty$ are listed in \ref{thm:BW_alg_length}.
\end{theorem}
In fact, the paper \cite{BW5} also classifies the Jacobi algebras with algebraic length that are infinite dimensional, and this gives bonus geometric applications to crepant divisor-to-curve contractions.

After the conclusion of the classification of flops, it turns out that  Koll\'ar's length is equal to the algebraic length, and so there is really only one invariant, defined in two completely different ways.  It remains convenient to list the classification as a table in which the rows are indexed by the length, as we do in \ref{thm:GeomReal}, partly since slightly different phenomena appear in the Type $E$ ($\ell\geq 3$) cases.

\subsection{Derived contraction algebras.}\label{subsec:Derived}
This section very briefly sketches the two approaches to proving the theoretical classification result in \ref{conj:DW}, and also its technical companion \ref{thm:DWThm}\eqref{thm:DWThm 2}. Throughout, let $\scrX\to\Spec\scrR$ be a smooth $3$-fold flopping contraction. 

\subsubsection{Unique enhancements.} The first approach, due to Jasso--Muro \cite{JKM}, proves \ref{conj:DW} via first establishing very strong properties of the singularity category $\scrT=\uCM\scrR$. In fact, their main result works much more generally: the following is a restricted form of the $d=2$ Derived Auslander--Iyama Correspondence, over the field $k=\mathbb{C}$.
\begin{theorem}[{\cite[Thm A]{JKM}}]\label{thm:JKM_thmA}
There are bijective correspondences between the following.
\begin{enumerate}
\item\label{thm:JKM_thmA 1} Equivalence classes of pairs $(\scrT,c)$ where
\begin{enumerate}
    \item $\scrT$ is a $\mathbb{C}$-linear, Hom-finite algebraic triangulated category in which idempotents split.
    \item $c\in\scrT$ is a basic $2\mathbb{Z}$-cluster-tilting object, namely $c[2]\cong c$ and furthermore
    \[
    \mathrm{add}\, c=\{x\in\scrT\mid \Hom_{\scrT}(x,c[1])=0\}=\{x\in\scrT\mid \Hom_{\scrT}(c,x[1])=0\}.
    \]
\end{enumerate}
\item\label{thm:JKM_thmA 2} Equivalence classes of pairs $(\Lambda,I)$ where
\begin{enumerate}
    \item $\Lambda$ is a basic self-injective finite dimensional $\mathbb{C}$-algebra, which is twisted periodic.
    \item $I$ is an invertible $\Lambda$-bimodule such that $\Omega^{4}_{\Lambda^{e}}\Lambda\cong I$ in the stable category of $\Lambda$-bimodules.
\end{enumerate}
\end{enumerate}
The bijection sends $(\scrT,c)\mapsto (\,\End_{\scrT}(c),\Hom_\scrT(c,c[-2])\,)$. Moreover, the algebraic triangulated categories in \eqref{thm:JKM_thmA 1} admit a unique \textnormal{DG}-enhancement.  
\end{theorem}
Now the base $\scrR$ of any smooth flopping contraction is isolated cDV, and thus an isolated hypersurface singularity. The first property implies that the Hom-spaces in $\scrT\colonequals \uCM\scrR$ are finite dimensional, and the second that $[2]\cong\mathrm{Id}$.  Further, by e.g.\ \ref{thm:AMCorr}, the contraction algebra $\Acon=\End_{\scrT}(c)$ for some cluster-tilting object $c$, which is automatically $2\mathbb{Z}$-cluster-tilting since $[2]\cong\mathrm{Id}$.

\medskip
To prove \ref{conj:DW}, as explained in \S\ref{subsec:DWtrue} we can reduce to the case where two contraction algebras are isomorphic. So, set $\scrT_1=\uCM\scrR_1$ and $\scrT_2=\uCM\scrR_2$, and suppose that $\Acon=\End_{\scrT_1}(c_1)\cong\End_{\scrT_2}(c_2)=\Bcon$.  Since $[2]\cong\mathrm{Id}$, the pairs $(\End_{\scrT_1}(c_1), \Hom_{\scrT_1}(c_1,c_1))$ and $(\End_{\scrT_2}(c_2), \Hom_{\scrT_2}(c_2,c_2))$ are equivalent, and so applying \ref{thm:JKM_thmA} yields a triangulated equivalence $\uCM\scrR_1\simeq\uCM\scrR_2$ which sends $c_1\mapsto c_2$.  By the final unique enhancement statement within \ref{thm:JKM_thmA}, this lifts to an isomorphism in the homotopy category of DG categories between the unique DG enhancements.  From here, by \cite{HuaKeller} the Tjurina algebras of both $\scrR_1$ and $\scrR_2$ can be obtained via taking Hochschild cohomology, and are isomorphic. Thus $\scrR_1\cong\scrR_2$ by Mather--Yau.

In the course of the proof of \ref{thm:JKM_thmA}, the $2\mathbb{Z}$-derived contraction algebra $\AconTwoZ\colonequals \mathbf{R}\mathrm{End}_{\scrT_1}(c_1)$, which is a DG-algebra concentrated in only even degrees, is shown to be determined by $\Acon$.  More precisely, \cite[Theorem 3.1.1]{JKM_overview} shows that two $2\mathbb{Z}$-derived contraction algebras whose zeroth cohomologies are isomorphic as algebras must necessarily be quasi-isomorphic. Since $\Aconder\cong\uptau_{\leq 0}(\AconTwoZ)$, \ref{thm:DWThm}\eqref{thm:DWThm 2} also follows.

\subsubsection{The return of the image.}
Alternatively, consider $\scrS=\bigoplus_{i=1}^n\scrO_{\Curve_i}(-1)\in\coh\scrX$, where $\Curve_1,\hdots,\Curve_n$ are the exceptional curves of $\scrX\to\Spec\scrR$.  The second approach to prove \ref{conj:DW} is to show that the DG-algebra $(\Aconder)^!=\mathrm{REnd}_{\scrX}(\scrS)$ can be directly reconstructed from the isomorphism class of $\mathrm{A}_{\mathrm{con}}$.  By Koszul duality \cite[Theorem 2.8]{KalckYang} and \ref{thm:Jenny}, this suffices.

We begin by recalling the solution to a strictly easier problem.  Given an injective ring homomorphism $\upvarphi\colon\mathrm{A}\hookrightarrow \mathrm{B}$, the first isomorphism theorem asserts that there is a third ring $\mathrm{Im}\upvarphi$, called the image, such that there is an isomorphism and commutative diagram as follows.
\[
\begin{array}{c}
\begin{tikzpicture}
\node (A) at (0,0) {$\mathrm{A}$};
\node (B) at (3,0) {$\mathrm{B}$};
\node (C) at (3,-1) {$\mathrm{Im}\upvarphi$};
\draw[->] (A)--(B);
\draw[->] (A)--node[below]{$\scriptstyle\cong$}(C);
\draw[->,densely dotted] (C)--(B);
\end{tikzpicture}
\end{array}
\]
Since $\mathrm{Im}\upvarphi$ is defined by all those elements that `come from'~$\mathrm{A}$, the statement $\mathrm{A}\cong\mathrm{Im}\upvarphi$ is really \emph{all} about $\mathrm{A}$. The ring $\mathrm{B}$ plays almost no role. This is not optimal if, say, you are trying to reconstruct $\mathrm{A}$ from $\mathrm{B}$.

Returning to the smooth $3$-fold flopping contraction $\scrX\to\Spec\scrR$, via the derived equivalence \eqref{VdBtilt}, we can replace $\scrS$ by its (finite) projective resolution $\scrP$ of $\Lambda\colonequals \End_{\scrR}(\scrR\oplus N)$-modules, and also consider the finite complex $\mathrm{P}=\scrP\otimes_{\Lambda}\Acon$ of $\Acon$-modules. There is a morphism of DG-algebras
\begin{equation}
\scrA:= \mathrm{REnd}_{\scrX}(\scrS)\stackrel{\scriptstyle \mathrm{qis}}{\simeq}
\mathrm{REnd}_{\Lambda}(\scrP)
\xrightarrow{\otimes_{\Lambda}\Acon} 
\mathrm{REnd}_{\Acon}(\mathrm{P})
=:\scrB\label{eqn:InjCohom}
\end{equation}
which (via some known moves) is injective on cohomology. It turns out that $\mathrm{A}_{\mathrm{con}}$ determines $\mathrm{P}$, it being the start of a minimal projective resolution of its simples, and thus $\mathrm{A}_{\mathrm{con}}$ builds $\scrB$.  The upshot is that there is a morphism of DG-algebras $\scrA\to\scrB$ which is injective on cohomology, and we are required  to reconstruct $\scrA$ from $\scrB$.  

Here the first isomorphism theorem analogy re-enters. In general,  reconstructing $\scrA$ from $\scrB$ should not be possible: $\scrB$ should be auxiliary information, and the statement should only concern $\scrA$.  However in the flops situation, if you believe \ref{conj:DW} to be true, then this strongly suggests that for the particular morphism \eqref{eqn:InjCohom} there must be an \emph{intrinsic} subalgebra of $\scrB$ (determined up to quasi-isomorphism of the DG-algebra $\scrB$), which acts as `the image' in this setting, from which $\scrA$ can be reconstructed.

This does exist, and since the morphism \eqref{eqn:InjCohom} is injective only on cohomology, to construct it we now pass to $A_\infty$-algebras. Replacing $\scrB$ by its strictly unital $A_\infty$-minimal model, we obtain the following.
\[
\begin{tikzpicture}
%\node (A) at (0,0) {$\scrB:=\mathrm{REnd}_{\mathrm{A}_{\mathrm{con}}}(\mathrm{P})$};
\node (A) at (2,0) {$\scrB$};
\node (B) at (5,0) {$(\mathrm{H}^*(\scrB),\{m_n\}_{n \geq 2})$};
\draw[->] ($(A)+(0.4,0.1)$) -- ($(B)+(-1.6,0.1)$)  ;
\draw[<-] ($(A)+(0.4,-0.1)$) -- ($(B)+(-1.6,-0.1)$)  ;
\end{tikzpicture}
\]
Inside the right-hand side is the subspace $\scrN_\scrB:=\mathbb{C}\mathsf{1}\oplus \mathrm{H}^{>0}(\scrB)$.
\begin{lemma}\label{lem:Easy}
With notation as above, the following statements hold.
\begin{enumerate}
\item\label{lem:Easy 1} For any\footnote{implicitly, we assume the very weak condition that $\mathrm{H}^*(\scrB)\neq 0$.} unital DG-algebra $\scrB$, the  $\mathrm{A}_\infty$-structure $\{m_n\}$ restricts to give an $\mathrm{A}_\infty$-structure on $\scrN_\scrB$.
\item\label{lem:Easy 2} Any quasi-isomorphism $\scrB_1\to\scrB_2$ of unital \textnormal{DG}-algebras restricts to an $\mathrm{A}_\infty$-isomorphism $\scrN_{\scrB_1}\to\scrN_{\scrB_2}$.
\end{enumerate}
\end{lemma}
With that, the proof that $\Acon$ determines $\scrA=\mathrm{REnd}_{\scrX}(\scrS)$ proceeds as follows \cite{KLW}. Given two flops, with contraction algebras $\mathrm{A}_{\mathrm{con}}$ and $\mathrm{B}_{\mathrm{con}}$ respectively, if $\Acon\cong\Bcon$ then
\[
\begin{tikzpicture}
\node (a) at (0,0) {$\scrB_1\colonequals \mathrm{REnd}_{\Acon}(\mathrm{P})$};
\node (b) at (4,0) {$\mathrm{REnd}_{\Bcon}(\mathrm{P'})=:\scrB_2$};
\draw[->] (a) --node[above]{$\scriptstyle \simeq$} (b);
\end{tikzpicture}
\]
and thus by \ref{lem:Easy}\eqref{lem:Easy 2} there is an $A_\infty$-isomorphism $\scrN_{\scrB_1}\to\scrN_{\scrB_2}$.  But since the morphisms \eqref{eqn:InjCohom} are injective on cohomology, $\scrN_{\scrB_1}$ is quasi-isomorphic to $\mathrm{REnd}_{\scrX_1}(\scrS)$, and $\scrN_{\scrB_2}$ is quasi-isomorphic to $\mathrm{REnd}_{\scrX_2}(\scrS')$. Combining these statements, $\mathrm{REnd}_{\scrX_1}(\scrS)$ and $\mathrm{REnd}_{\scrX_2}(\scrS')$ are quasi-isomorphic, which is sufficient to prove \ref{conj:DW}. Further, since by Koszul duality $\Aconder\simeq\mathrm{REnd}_{\scrX_1}(\scrS)^!$, and likewise $\Bconder\simeq\mathrm{REnd}_{\scrX_2}(\scrS')^!$, it also follows that $\Aconder$ is quasi-isomorphic to $\Bconder$, and thus \ref{thm:DWThm}\eqref{thm:DWThm 2} also follows.

\section*{Acknowledgments.}
Many thanks are due to all my co-authors, colleagues, postdocs, and students who have supported, taught and inspired me, particularly Osamu Iyama, Will Donovan and Gavin Brown for their generous collaborations on the topic of this survey for over a decade.  Thanks also to Aidan Schofield, who in the mists of time both began and actively encouraged noncommutative dreaming.

The author is funded by ERC Consolidator Grant 101001227 (MMiMMa), and previously by EPSRC grants EP/R034826/1, EP/R009325/1 and EP/K021400/2.

\end{document}